\magnification\magstep1
\baselineskip=18pt
\def\w{\widetilde}
\def\i{\infty}
\def\v{\vert}
\def\V{\Vert}
\centerline{{\bf  Interpolation Between
$H^p$ Spaces and Non-Commutative
Generalizations I}\footnote*{Supported in part by N.S.F.
grant DMS 9003550}} \vskip12pt \centerline {by  Gilles
Pisier} \vskip12pt

{\bf Abstract}  
We give an elementary proof  that the $H^p$ spaces over
the unit disc (or the upper half plane) are the 
interpolation spaces for the real method of interpolation
between $H^1$ and $H^{\infty}$. This was originally
proved by Peter Jones. The proof uses only the boundedness
of the Hilbert transform and the classical factorisation
of a function in $H^p$  as a product of two functions in
$H^q$ and $H^r$  with ${1/q}+{1/r}=1/p$. This proof
extends without any real extra difficulty to the
non-commutative setting and to several Banach space
valued extensions of $H^p$ spaces. In particular, this
proof easily extends  to the 
couple $H^{p_{0}}(\ell_{q_0}),H^{p_{1}}(\ell_{q_1}) $,
with $1\leq p_{0}, p_{1}, q_{0}, q_{1} \leq \infty$. In
that situation, we prove that 
 the real
interpolation spaces and the K-functional are induced ( up
to equivalence of norms ) by the same objects for the
couple  $L_{p_0}(\ell_{q_0}),L_{p_1}(\ell_{q_1}) $. In an
other direction,   let us denote by $C_p$  the space of all
compact operators $x$ on
 Hilbert space such that $tr(|x|^p) <\infty$. Let $T_p$
be the subspace of all upper triangular matrices
relative to the canonical basis.
If
$p=\infty$, $C_p$ is just the space of all compact
operators. Our proof allows us to show for instance that 
the space $H^p(C_p)$ (resp. $T_p$) is the interpolation
space of parameter $(1/p,p)$ between $H^1(C_1)$
(resp. $T_1$) and $H^\infty(C_\infty)$ (resp. $T_\i$). We
also prove a similar result for the complex 
interpolation method. Moreover,  extending a recent
result of Kaftal-Larson and Weiss, we prove that the
distance to the subspace of upper triangular matrices in
$C_1$ and $C_\infty$ can be essentially realized
simultaneously by the same element.

\vfill\eject
\noindent {\bf  Introduction.}

\magnification\magstep1
\baselineskip=18pt
\magnification\magstep1
\baselineskip=18pt

\def\i{\infty}
\def\l{\ell}
\def\t{\theta}
\def\la{\lambda}
\def\ra{\rightarrow}
\def\i{\infty}

\def\n{\noindent}

Let $0<p\leq\i$. We will denote simply by $L_p$ the
$L_p$-space relative to the circle group ${\bf T}$
equipped with its normalised Haar measure denoted by $m$.
We will denote by $H^p$ the classical Hardy space of
analytic functions (on the unit disc $D$ of the complex
plane). It is well known that this space can be identified
with a closed subspace of  $L_p$, namely the closure in
$L_p$ (for $p=\i$ we must take the weak*-closure) of the
linear span of the functions $\{e^{int}|n\ge 0\}$. We
refer e.g. to [G] or [GR] for more information on
$H^p$-spaces.

Let us recall here the definitions of the $K_t$ and $J_t$
functionals which are fundamental in the real
interpolation method.
Let $A_0,A_1$ be a compatible couple of Banach (or
quasi-Banach) spaces. This just means that $A_0,A_1$ are
continuously included into a larger topological vector
space (most of the time left implicit), so that we can
consider unambiguously the sets $A_0+A_1$ and $A_0\cap
A_1$. 

\n For all $x\in A_0+A_1$ and
for all $t>0$, we let $$K_t(x,A_0,A_1)= \inf
\big({\|x_0\|_{A_0}+t\|x_1\|_{A_1}\ | \ x=x_0+x_1,x_0\in
A_0,x_1\in A_1}).$$ For all $x\in A_0\cap A_1$ and for
all $t>0$, we let
$$J_t(x,A_0,A_1)=\max(\|x_0\|_{A_0},t\|x_1\|_{A_1}).$$
Recall that the (real interpolation) space
$(A_0,A_1)_{\theta,p}$ is defined as the space of all $x$
in $A_0+A_1$ such that $\|x\|_{\theta,p} <\infty$ where
 $$\|x\|_{\theta,p} =(\int{(t^{-\theta}K_t(x,A_0,A_1))^{p}
dt/t})^{1/p} .$$
We also recall that there is a parallel definition of $(A_0,A_1)_{\theta,p}$
using the $J_t$ functional which leads to the same Banach
space with an equivalent norm.

\n For example, if $1\le {p_0},{p_1}, q \le \i$ and
$0<\t<1$, we have
$$ (L_{p_0},L_{p_1})_{\t,q} = L_{p,q}$$
where $ L_{p,q}$ is the classical Lorentz space,
identical to  $L_{p}$
if $p=q$. We refer to [BL] for more details.

In section 1, we give a new proof of the following
interpolation theorem of Peter Jones [J1], as reformulated
by Sharpley (cf. [BS] p. 414):

\n There is a constant $C$ such that
$$ \forall f \in H^1+H^\i, \quad\forall t>0,\quad
K_t(f,H^1,H^\i) \leq C
K_t(f,L^1,L^\i).$$

\n We should recall that 
$$K_t(f,L^1,L^\i)=\int_0^t{f^{*} ds}=\sup\{\int_E |f| dm
|\ \  E\subset {\bf T}\ , \  m(E)=t \} ,$$ where we have
denoted by $f^{*}$ the non-increasing rearrangement of the
function $|f|$.  The difficulty of Jones' theorem lies in
the fact that the optimal decomposition which realizes
$K_t(f,L^1,L^\i)$ is obtained by truncating the function
$f$. If $f$ is analytic, this operation clearly spoils the
analyticity, and the problem is to find a substitute,
something like a truncation but which preserves
analyticity.

We should mention that a different proof of
Jones'results (including some results which cannot be
obtained by our method) has already been obtained  a few
months ago by Quanhua Xu. However, Xu's argument does
not seem to extend to the non-commutative case.

We now describe the (very simple) method of proof we use
throughout this paper, we call it the "square/dual/square"
argument.

\n Let us say that the couple  $(H^p,H^q)   $ is K-closed
if there is a constant $C$ such that
$$ \forall f \in H^p+H^q, \quad\forall t>0,\quad
K_t(f,H^p,H^q) \leq C
K_t(f,L^p,L^q).$$

\n The first step consists in showing the following
"squaring" property:

\n (0.1) \quad If $(H^{2p},H^{2q}) $ is K-closed
and if the pointwise product defines a bounded bilinear map
from $H^{2p}\times H^{2q}$ into the interpolation space 
$(H^p,H^q)_{1/2,\i}$, then $(H^p,H^q)$ is K-closed.

\n The next step is a dualisation (this seems to be the
point that has been overlooked by previous researchers)

\n (0.2) \quad  The couple $(H^p,H^q)$ is K-closed iff
$(H^{p'},H^{q'})$ is also K-closed ($1\le p,
q\leq\i,1/p+1/p'=1/q+1/q'=1).$

\n We can then sketch our "square/dual/square"-proof of
the fact that $(H^1,H^\i)$ is K-closed as follows:

\n By (0.1), it suffices to show that $(H^2,H^\i)$ is
K-closed, then by (0.2) it suffices to show that
$(H^2,H^1)$ is K-closed, but then by (0.1) again, it
suffices to show that $(H^4,H^2)$ is K-closed, and this is
an obvious and well known consequence of Marcel
Riesz' theorem on the simultaneous boundedness of the
Hilbert transform on  $L_p$ for all $1<p<\i$.

Our proof emphasizes the existence of a "simultaneous good
approximation" to $H^1$ and $H^\i$. More precisely, we have

\n (0.3) There is a constant $C$, such that for all $f\in
L_\i$, there is a function $h\in H^\i$ such that we have
simultaneously $$\|f-h\|_\i \leq C\  dist_{L_\i} (f,H^\i)  
\quad   and \quad \|f-h\|_1 \leq C\  dist_{L_1} (f,H^1) .$$
 
As far as we know at the time of this writing, these
results are known only in dimension 1, and are open in
higher dimension either for the ball or the polydisc. We
refer the reader to  [J2] for a survey of what is known
in the latter case.

In
section 2, we prove a non-commutative analogue of Peter
Jones' theorem, where the space $L^p$ is replaced by the
space $C_p$ of all compact operators $x$ on $\l_2$ such
that $ tr|x|^p <\i $, and $H^p$ is replaced by the
subspace $T_p$ of all upper triangular matrices. 

This result, which was motivated by and which improves a
result of [KLW], says again that the $K_t$-functional for
the couple $(T_1,T_\i)$ is induced (up to a constant
independent of $t$) by the $K_t$-functional for
the couple $(C_1,C_\i)$. As a corollary, we identify the
real interpolation spaces for the couple $(T_1,T_\i)$.

In that case also, there is a simultaneous good
approximation to $T_1$ and $T_\i$, as in (0.3) above.

In section 3, we discuss the case of Banach space valued
$H^p$-spaces. In particular, we show that Jones' theorem
is also true for the couple of operator valued $H^p$-spaces
$(H^1(C_1), H^\i(C_\i))$, i.e. the $K_t$-functional is
induced (up to a constant
independent of $t$) by the $K_t$-functional for
the couple $(L_1(C_1),L_\i(C_\i))$. (See
Theorem 3.3 and its corollaries for more precision.) As a
consequence, we
can again identify the real interpolation spaces.  This
result is closely related to the result in section 2. (In
fact, one can deduce from it the  above result on
$(T_1,T_\i)$.) More generally, we obtain similar results
for the couples $(H^{p_0}(\ell_{q_0}),H^{p_1}(\ell_{q_1})
)$, and
 $(H^{p_0}(C_{q_0}),H^{p_1}(C_{q_1})) $,
with $1\leq p_{0}, p_{1},q_{0},q_{1} \leq \infty$. We
doubt
that Jones' proof can be adapted to all these
cases. In the case of
$(H^\i(\l_2),H^\i(\l_\i))$, our argument leads to a new
proof of a theorem of Bourgain [B], but we chose to write
this separately, we refer to [P].

In section 4, we consider similar problems for the
complex interpolation method.  Peter Jones [J1] proved
that $H^p=(H^1,H^\i)_\t$, with $1/p =1-\t$. Using what
seems to be a new idea in this context, we show that this
result can be deduced from a slightly extended version of
the real case. Our argument  extends to the
non-commutative case and gives $T_p=(T_1,T_\i)_\t$.

Although we state and prove our results on the unit
disc, there is no problem to extend them to the case of
the upper half plane. We leave this to the reader.

We now introduce a specific notation needed  to treat the
Banach space valued case. Let {\bf T} be the circle group
equipped with its normalized Haar measure $m$. Let
$0<p\le \i$. When $B$ is a Banach space, we denote by
$L_p(B)$ the  usual space of Bochner-$p$-integrable
$B$-valued functions on $({\bf  T},m)$, so that when
$p<\i$, $L_p \otimes B$ is dense in $L_p(B)$. We denote
by $\w{H}^p(B)$ the closure in $L_p(B)$ of all the finite
sums of the form $\sum_{0<k<n} x_k e^{ikt}$ with $x_k \in
B$. In other words, if we denote by ${\cal T}$ the space
of all analytic trigonometric polynomials, $\w{H}^p(B)$
is the closure in $L_p(B)$ of ${\cal T}\otimes B$. We
reserve the notation $H^p(B)$ (and simply $H^p$ if $B$ is
one dimensional) for the Hardy space of B-valued analytic
functions $f$ such that 
$$\sup_{r<1}(\int\|f(re^{it})\|^p dm(t))^{1/p} <\i.$$
Again, see [G,GR] for more information.

When $B$  is reflexive,   is  a separable
dual or is an $L_1$-space (in particular  if $B$ is  
finite dimensional), then
it is well known that $\w{H}^p(B)=H^p(B)$ for all 
$p<\i$,  and ${H}^p(B)$ can be identified with a subspace
of $L_p(B)$ for all $p\le\i .$  We refer to
[BuD,E,HP] for more information on this property, called
the analytic Radon-Nikodym property.

The next proposition although very simple will be
essential in the sequel. We suspect that the importance of
the equivalence $(i)\Leftrightarrow (ii)$  has been
overlooked although its proof is routine. We should
emphasize that the exponents $p,q$ in $(i)$ and $(ii)$ are
the same, they are $\underline{not}$ conjugate to each
other.  

\proclaim  Proposition 0.1. Let
$1\leq p\leq q\leq \infty$. Consider an interpolation
couple
of Banach spaces $(A_0,A_1)$,  the following
are equivalent: 
$(i)$ There is a constant $C'$ such that
$$ \forall f \in \w{H}^p(A_0)+\w{H}^q(A_1), \quad\forall
t>0,\quad K_t(f,\w{H}^p(A_0),\w{H}^q(A_1)) \leq C'
K_t(f,L^p(A_0),L^q(A_1)).$$
$(ii)$ There is a constant $C$ such that
$$ \forall f \in {[L^p(A_0)/ \w{H}^p(A_0)]}\cap
{[L^q(A_1)/\w{H}^q(A_1)]},\quad \forall t>0,\quad  \exists
\hat{f} \in L^p(A_0)\cap L^q(A_1)$$
 satisfying
$$J_t(\hat{f},L^p(A_0),L^q(A_1))
\leq C
J_t(f,L^p(A_0)/{\w{H}^p(A_0)},L^q(A_1)/{\w{H}^q(A_1)}).$$
$(iii)$ There is a constant $C$ such that $$ \forall f \in
{[L^p(A_0)/ \w{H}^p(A_0)]}\cap {[L^q(A_1)/\w{H}^q(A_1)]}, 
\quad  \exists \hat{f} \in L^p(A_0)\cap L^q(A_1)$$
 satisfying
$$\|\hat{f}\|_{L^p(A_0)}\leq C\|f\|_{L^p(A_0)/\w{H}^p(A_0)},$$
$$\|\hat{f}\|_{L^q(A_1)}\leq C\|f\|_{L^q(A_1)/\w{H}^q(A_1)}.$$

\medskip
In the above statement we regard the spaces 
$L^p (A_0)/\w{H}^p (A_0)$ and $L^q(A_1)/\w{H}^q(A_1)$ as included
via the Fourier transform $f\ra
(\hat{f}(-1),\hat{f}(-2),\hat{f}(-3),...)$ in the space of
all sequences in $A_0 + A_1$. In this way, we may view
these quotient spaces as forming a compatible couple for
interpolation. (For the subspaces 
$\w{H}^p(A_0),\w{H}^q(A_1)$, there is no problem, we may clearly
consider them as a compatible couple in the obvious
way.) 
\medskip
\n {\bf Proof:} For brievity, we will denote simply
$L^p/{H}^p(A_0)$ instead of $L^p(A_0)/\w{H}^p(A_0)$, we
will also write $L^p$,$H^p$,..instead of $L^p(A_0)$ ,
$\w{H}^p(A_0)$...no confusion should arise.
The proof is routine. We indicate first the argument for
$(i) \Rightarrow (ii)$ which is the one we use below.

Assume $(i)$.
Let $f$ be as above such that
 $J_t(f,L^p/H^p(A_0),L^q/H^q(A_1))<1$. Then let $g_p \in
L^p(A_0)$\quad  and $g_q \in L^q(A_1)$ be such that 
$$\|g_p\|_{L^p} <1, \quad \|g_q\|_{L^q} <t^{-1}, \quad 
f=g_p+H^p(A_0) ,\quad f=g_q+H^q(A_1).$$
 Therefore, $g_p-g_q$ must be in $H^p+H^q$ and 
 $$K_t(g_p-g_q,L^p(A_0),L^q(A_1)) \leq
\|g_p\|_{L^p}+t\|g_q\|_{L^q}<2.$$
By $(i)$, we have $K_t(g_p-g_q,H^p,H^q)<2C'$, hence there
are $h_p\in H^p(A_0)$ and $h_q\in H^q(A_1)$ such that
$g_p-g_q=h_p-h_q$ and $\|h_p\|_{H^p}+t\|h_q\|_{H^q}<2C'$.
Now if we let $\hat f =g_p-h_p=g_q-h_q$,
then we find that $\hat{f} \in L^p(A_0)\cap
L^q(A_1)$,$f=\hat f +H^p(A_0)$ in the space $L^p/H^p(A_0)$ and 
$f=\hat f +H^q(A_1)$ in the space $L^q/H^q(A_1)$ and
moreover $$J_t(\hat f,L^p,L^q) \leq \max(\|\hat
f\|_{L^p},t\|\hat f\|_{L^q})\leq 1+2C'.$$ By homogeneity
this completes the proof of $(i)\Rightarrow (ii)$ with
$C\leq 1+2C'$. The converse is similar, we skip the
details.
The implication $(ii) \Rightarrow (iii)$ is easy, just take
$$t=(\|f\|_{L^p(A_0)/H^p(A_0)})(\|f\|_{L^q(A_1)/H^q(A_1)})^{-1}.$$
The converse $(iii) \Rightarrow (ii)$ is trivial. q.e.d.

\medskip

{\bf Remark:}The preceding statement would also remain
valid if we had defined $\w{H}^{\i}(B)$ as the subspace
of $L_{\i}(B)$ formed by the functions with a Fourier
transform vanishing on the negative integers. See the
end of this section for a more general viewpoint.
 \medskip

We recall the following basic fact:
If
$1<p_0<p_1<\infty$ then there is a constant $C$  such that
for all $t>0$ we have:
$$  \quad    \forall f \in H^{p_0}
+H^{p_1},\quad \forall t>0,\quad K_t(f,H^{p_0} ,H^{p_1})
\leq  CK_t(f,L^{p_0} ,L^{p_1}) .\leqno(0.4) $$

This is an obvious consequence of the simultaneous 
boundedness of the orthogonal projection $P:L^2\rightarrow
H^2$ on all the $L^p$ spaces (or equivalently of the same
for the Hilbert transform).
This "simultaneous" boundedness of $P$ obviously also
implies that if \ \ \ \ \ \ \qquad $1<p_0<p<p_1<\infty$ and
if $1/p=(1-\theta)/p_0 + \theta/p_1 $, we have
$$H^p =(H^{p_0} ,H^{p_1})_{\theta,p}$$ or more
generally, if we define $H^{p,q}$ as the space of
analytic functions in the disc with boundary values in
the Lorentz space $L^{p,q}$ (on the circle), then, for any
$1\leq q\leq \i$,  we have
$$H^{p,q} =(H^{p_0}
,H^{p_1})_{\theta,q},$$
and in particular
$$H^p \subset H^{p,\i} = (H^{p_0}
,H^{p_1})_{\theta,\i}.\leqno(0.5)$$

\vfill\eject

 \n{\bf  1. The proof of Peter Jones' Theorem.}

We prove the theorem in several steps: starting from (0.4)
restricted to $p_0,p_1$ both finite and more than $1$, 
we will progressively extend the set of couples
$(p_0,p_1)$ for which (0.4) is valid until we eventually
have eliminated all restrictions on $p_0,p_1$.

\proclaim Proposition 1.1. For all $1<p<q<\infty$
 we have
$$ \quad H^p\subset (H^1,H^q)_{\theta,\infty},\leqno(1.1)$$
with norm bounded by some constant $K(p,q)$, where $
0<\theta<1$ satisfies  $1/p=1-\theta+\theta/q.$ 

\medskip
\n {\bf Proof:} Choose any number $r>q$, and define $r'$, $s$
and $t$ by the relations
$$ {1/r}+{1/r'}=1, \quad {1/r}+{1/s}=1/p,\quad
{1/r}+{1/t}=1/q .$$
Observe that ${1/s}=(1-\theta)/{r'}+\theta/t.$ Let f be in
the unit ball of $H^p$, write $f=gh$ with $g$ and $h$
respectively in the unit balls of $H^r$ and $H^s$. By the
above basic fact (0.4) we have
$$H^s \subset (H^{r'} ,H^{t})_{\theta,\infty}$$ 
and this
inclusion has norm less than (say) $C$. 
Observe that the operation of multiplication by $g$
maps (by H\"older's inequality) the unit ball of $H^{r'}$
(resp. $H^{t}$) into that of $H^1$
(resp. $H^q$), hence it maps the unit ball of $(H^{r'}
,H^{t})_{\theta,\infty}$ into that of
$(H^1,H^q)_{\theta,\infty}$.
Therefore the norm of $f=gh$ in
the space $(H^1,H^q)_{\theta,\infty}$ is less than $C$,
which completes the proof. (This statement is also
immediate using the complex interpolation method.) q.e.d.
\medskip

 The proof of the next proposition, although very simple
is important in the sequel.
\medskip
\proclaim  Proposition 1.2. For each
$1<q<\infty$, there is a constant $C'$ such that $$\forall
t>0,\quad \forall f \in H^1+H^q, \qquad K_t(f,H^1,H^q) \leq
C'K_t(f,L^1,L^q).$$\medskip

\n {\bf Proof:}  Let $f$ be analytic in the disc and such
that $K_t(f,L^1,L^q) <1$.We factorize $f$ as $f=BF^{2}$
with $F$ non vanishing and $B$ a Blaschke product. Then 
since $|F|= |f|^{1/2}$ on the unit circle, we clearly have
$K_{t^{1/2}}(F,L^2,L^{2q}) <2^{1/2}$, hence by (0.1), 
$K_{t^{1/2}}(F,H^2,H^{2q})<2^{1/2}C$. Therefore, there are
$\underline{analytic}$ functions $g_0$ and $g_1$
such that  $$  \quad F=g_0+g_1 \qquad
\|g_0\|_2+t^{1/2}\|g_1\|_{2q} <2^{1/2}C.\leqno(1.2)$$
Now we can write $f=B(g_0+g_1)^2=B(g_0^2+g_1^2+2g_0g_1)$,
hence$$ \quad K_t(f,H^1,H^q) \leq
K_t(g_0^2+g_1^2,H^1,H^q) + K_t(2g_0g_1,H^1,H^q).\leqno(1.3)$$
By (1.3) we have 
$$ \quad K_t(g_0^2+g_1^2,H^1,H^q)\leq 2C^2,\leqno(1.4)$$ 
and on the other hand by Holder $ \|2g_0g_1\|_p \leq
2C^2t^{-1/2}$ where $1/p=1/2+1/{2q}$.
Note that $1<p<q$, and that
$1/p=1-\theta+\theta/q$ with $\theta =1/2$ , so that by
Proposition 1.1  for some constant $K$ we have
$$\|2g_0g_1\|_{(H^1,H^q)_{\theta,\infty}} \leq
K2C^2t^{-1/2}. $$
Hence,in particular,
 $t^{-\theta}K_t( 2g_0g_1,H^1,H^q) \leq K2C^2t^{-1/2}$,
so that  $$  \quad K_t( 2g_0g_1,H^1,H^q) \leq
K2C^2.\leqno(1.5)$$
Returning to (1.3),we see that (1.4) and (1.5) imply
$$ K_t(f,H^1,H^q)\leq 2C^2+K2C^2.$$ q.e.d.\medskip

\n{\bf Remark:} At this point, we can easily check (0.1)
by a minor modification of the preceding proof. We will
 refer to (0.1) in the 
sequel as "the
squaring argument".

The
special nature of the K and J functionals on one hand and
of $H^p$ and $(H^p)^\perp$ on the other hand imply that
Proposition 1.2 has the following consequence.

\proclaim  Proposition 1.3. For each $1<q<\infty$,
there is a constant $C_q$ such that $$\forall
t>0,\quad \forall f \in {L^1/H^1}\cap{L^q/H^q},$$
$$\exists \hat{f} \in L^1\cap L^q  \quad satisfying\quad
 J_t(\hat f,L^1,L^q)\leq C_{q} J_t( f,L^1/H^1,L^q/H^q).$$
Equivalently,  $\forall f \in {L^1/H^1}\cap{L^q/H^q},$
$$ \exists \hat{f} \in L^1\cap L^q  \quad satisfying\quad
\|\hat f\|_{L^1}\leq C_q \|f\|_{L^1/H^1},\quad \|\hat
f\|_{L^q}\leq C_q \|f\|_{L^q/H^q}.\leqno(1.6)$$\medskip

\n {\bf Proof:} By proposition 0.1, this follows from
Proposition 1.2. q.e.d. \medskip

Up to now we have not used the duality between the $K_t$
and $J_t$ functionals, we now do so. We record below the
dual versions of the preceding two propositions.

\proclaim  Proposition 1.2*. For each
$1<p<\infty$, there is a constant $C'_p$ such that
\medskip$\forall t>0,\quad\forall f \in L^\infty /H^\infty
\cap L^p/H^p,$
 $$  \exists \hat f \in L^\infty \cap L^p \quad
satisfying\quad  J_t(\hat f,L^\infty ,L^p )
\leq C'_p J_t(f,L^\infty /H^\infty ,L^p/H^p).$$\medskip

\proclaim Proposition 1.3*. For each
$1<p <\infty$, there is a constant $C_p$ such that
$$\forall t>0,\quad \forall f \in H^\infty +H^p, \qquad
K_t(f,H^{\infty},H^p) \leq C_{p}K_t(f,L^\infty,L^p).$$
\medskip

\qquad The  proof is obvious, we just recall that if $p$
and $q$ are conjugate then the dual of the space
$K_{t}(L^1/H^1,L^p/H^p) $   
(resp.$ J_{t}(L^1/H^1,L^p/H^p)$) is isometrically
identifyable with the space\qquad $J_{t}(H^{\infty},H^q)$
  (resp.$K_{t}(H^{\infty},H^q)$), and that an
injection is an isomorphic embedding iff its adjoint is
onto (the relevant constants being the same).

Let us record here an immediate consequence of Proposition
1.3 and Proposition 1.2*

\proclaim Corollary1.4. Assume  $1<p<q<\i$
with $1/p={(1-\beta)}/1+\beta/q$ and
with $1/q={(1-\gamma)}/\infty+\gamma/p$. Then for all
$f\in {L^1/H^1}\cap{L^q/H^q}$   $$\|f\|_{L^p/H^p} \leq C_q
(\|f\|_{L^1/H^1})^{1-\beta}
.(\|f\|_{L^q/H^q})^{\beta} .\leqno(1.7)$$
 And for all $f\in
{L^p/H^p}\cap{L^\infty/H^\infty}$  
$$\|f\|_{L^q/H^q} \leq C'_{p}
(\|f\|_{L^p/H^p})^{\gamma}
.(\|f\|_{L^\infty/H^\infty})^{1-\gamma} .\leqno(1.8)$$
\medskip

\n {\bf Proof:} (1.7) follows immediately from (1.6),
and (1.8) can be proved similarly using Proposition 1.2*
instead of Proposition 1.3. q.e.d.

We now use a  simple "extrapolation" trick to obtain

\proclaim Proposition 1.5. For each
$1<p<\i$, there is a constant $K_p$ such that  $$\forall f
\in L^\infty /H^\infty \cap L^1/H^1  \qquad
 \|f\|_{L^p/H^p}
\leq K_p(\|f\|_{L^1/H^1})^{1-\theta}
.(\|f\|_{L^\infty/H^\infty})^{\theta},$$ where
$1/p=(1-\theta)/1+\theta/\i.$\medskip

\n {\bf Proof:} Combining (1.7) and (1.8), we find, 
$$\|f\|_{L^p/H^p} \leq C_q
(\|f\|_{L^1/H^1})^{1-\beta}
.(C'_{p}
{\|f\|_{L^p/H^p}}^{\gamma}
.{\|f\|_{L^\infty/H^\infty}}^{1-\gamma})^{\beta}.  $$
Hence, 
$${\|f\|_{L^p/H^p}}^{1-\beta\gamma} \leq C_q
(C'_{p})^\beta {\|f\|_{L^1/H^1}}^{1-\beta}
{\|f\|_{L^\infty/H^\infty}}^{(1-\gamma)\beta}.$$
This yields the desired inequality with
$$K_p=(C_q
(C'_{p})^\beta)^{1/{1-\beta\gamma}},$$
since
$1-\theta=(1-\beta)(1-\beta\gamma)^{-1}.$ q.e.d.\medskip

We can now complete our proof of Peter Jones' theorem.
(A different proof has already been given a few months
ago by Quanhua Xu [X3].)

\proclaim Theorem 1.6. There is a constant
$C$  such that for all $t>0$ we have:
$$    \forall f \in H^{1} +H^{\infty},\quad
\forall t>0,\quad K_t(f,H^{1} ,H^{\infty}) \leq 
CK_t(f,L^{1} ,L^{\infty}) . $$\medskip

\n {\bf Proof:} We simply reproduce the proof of Proposition
1.2, but this time we can take $q=\infty$ because of
Proposition 1.3* (applied to the case $p=2$). Moreover,
Proposition 1.5 allows us to complete that same proof
because by duality Proposition 1.5
 is equivalent to the assertion  $$H^{p'} \subset
(H^1,H^{\infty})_{1-\theta,\infty},$$
with norm bounded by some constant $K$. In other words,
Proposition 1.1 remains valid for $q=\i$. It is then easy
to complete the proof by the squaring argument of
Proposition 1.2. q.e.d.\medskip
\magnification\magstep1
\baselineskip = 18pt
\def\n{\noindent}

\proclaim Corollary 1.7.
For all $0 < \theta < 1, 1 \le q \le \infty$ we have
$$H^{pq} = (H^1, H^\infty)_{\theta q}$$
\n where ${1\over p} = {1-\theta\over 1} + {\theta \over \infty}$.
\medskip

\n {\bf Remark 1.8:} \ To prove Proposition 1.5, we can
invoke Tom Wolff's interpolation theorem [W]. Indeed,
Proposition~1.3 gives us the ``right answer'' for the
interpolation spaces $(L_1/H^1, L_q/H^q)$ for $q<\infty$
and Proposition~1.2  gives us the case $(L_p/H^p,
L_\infty/H^\infty)$ for $p>1$. Actually, Corollary~1.7 can
be deduced directly from Propositions~1.2 and 1.3$^*$
using Wolff's results in [W].\medskip

\n {\bf Remark 1.9:} \ It is easy to extend Theorem~1.6 to the case of $H^r$
with $0 < r < 1$. First, we can check

$$H^1 \subset (H^r, H^\infty)_{\alpha,\infty}\leqno (1.9)$$

\n with ${1\over 1} = {1-\alpha\over r} + {\alpha \over \infty}$. Indeed,
we easily prove $H^1\subset (H^r, H^q)_{\alpha,\infty}$ with ${1\over 1} =
{1-\alpha\over r} + {\alpha \over q}$ for $q<\infty$ by the same method as
above. Then we can obtain (1.9) from Wolff's theorem [W]. Using (1.9), it
is immediate to adapt the preceding arguments to prove 
Theorem~1.6 with
$H^r$ and $L^r (0 < r < 1)$ instead of $H^1$ and $L^1$.

\n {\bf Remark 1.10:} (i) The same method will prove that
the couple of quasi-normed spaces $(H^{1,\i},H^\i)$ is
K-closed relative to $(L^{1,\i},L^\i)$. The same argument
works. Note however that we already know a priori from
Corollary 1.7 that the real interpolation spaces between
$(H^{1,\i},H^\i)$ are  the same (by reiteration) as those
between $(H^{1},H^\i)$. Indeed, the inclusion
$(H^{1,\i},H^\i)_{\t,q} \subset H^{p,q}$ is the trivial
direction, and Corollary 1.7 provides the converse. A
fortiori the same is true for the interpolation spaces
between $(L_1/\bar{H}^{1}_0 ,H^\i)$.
(ii) By Holmstedt's formula (cf.[BL], p.52-53), it follows
from Jones' theorem that all the couples $(H^p,H^q)$ are
K-closed, for any $0<p,q\le \i$, and similarly for
couples of Lorentz spaces $(H^{p_0,q_0},H^{p_1,q_1})$ with
$p_0 \ne p_1$ . \medskip

Let us recapitulate
and at the same time formalize the preceding argument.

\n Consider a compatible couple $(A_0, A_1)$ of Banach
(or quasi-Banach) spaces. Assume given a closed subspace
$S\subset A_0 + A_1$ and let

$$S_0 = S\cap A_0,\qquad S_1  = S\cap A_1.$$

\n Let $Q_0 = A_0/S_0$ and $Q_1 = A_1/S_1$ be the associated
quotient spaces. Clearly $(Q_0, Q_1)$ form a
compatible couple since there are natural inclusion
maps $$Q_0 \to (A_0 +A_1)/S \quad {\rm and} \quad Q_1 \to
(A_0 + A_1)/S.$$
 \n We will say that the couple $(S_0, S_1)$ is
$K$-closed (relative to $(A_0, A_1)$) if there is a
constant $C$ such that
$$\eqalign{&\forall t>0 \quad \forall x\in S_0 + S_1\cr
&K_t(x, S_0, S_1) \le C K_t(x, A_0, A_1).}$$

\n We will say that $(Q_0, Q_1)$ is $J$-closed if for some constant $C$ we
have $$\forall t >0 \quad \forall x\in Q_0\cap Q_1 \quad
\exists \hat x \in A_0 \cap A_1 \ \ {\hbox{such that}}\ \ 
J_t(\hat x, A_0, A_1) \le  J_t(x, Q_0, Q_1).$$

\n By the same argument as in Proposition 0.1 above one can show
that this is equivalent to the following `` simultaneous
lifting property'': \ $\forall x\in Q_0\cap Q_1 \quad
\exists \hat x \in A_0\cap A_1$ such that

$$x = \hat x + S_0 \quad {\rm in}\quad Q_0, \quad x = \hat x + S_1 \quad
{\rm in}\quad Q_1$$

\n and $$\|\hat x\|_{A_0} \le C\|x\|_{Q_0}, \quad \|\hat
x\|_{A_1} \le C \|x\|_{Q_1}.$$ 
Our terminology is motivated by the fact that, roughly
speaking,  $(S_0, S_1)$ is $K$-closed iff $S_0+ S_1$ is
closed in $A_0+A_1$ with a uniformity over $t$, while 
$(Q_0, Q_1)$ is $J$-closed iff $Q_0\cap  Q_1$ is closed in
$(A_0+A_1)/S$ with a uniformity over $t$. Then our key
observation in the preceding proof can be reformulated more
``abstractly'' as follows:\medskip

\proclaim  Proposition 1.11. $(S_0, S_1)$ is $K$-closed iff
$(Q_0, Q_1)$ is $J$-closed. \medskip

We leave the routine proof to the reader.

\n {\bf Remark 1.12:} \ Let us denote
$A_{\t,p}=(A_0,A_1)_{\t,p}$ and
$S_{\t,p}=(S_0,S_1)_{\t,p}$ . Assume that $(S_0, S_1)$ is
$K$-closed (relative to $(A_0,A_1)$). Then $S_{\t,p}$
can obviously be identified with a subspace of
$A_{\t,p}$  and  the  norm induced by $A_{\t,p}$ on 
$S_{\t,p}$ is
equivalent to the norm  of $S_{\t,p}$.
Moreover, the Holmstedt reiteration formula (cf. [BL],
p.52-53) for the K-functional shows that if
$0<\t_0\ne\t_1<1$, and if $1\le p_0, p_1 \le\i$, then the
couple $(S_{\t_0,p_0},S_{\t_1,p_1})$ is a fortiori
K-closed relative to $(A_{\t_0,p_0},A_{\t_1,p_1})$, and
also the couples $(S_0,S_{\t_1,p_1})$ and
$(S_{\t_0,p_0},S_1)$ are K-closed relative to respectively
 $(A_0,A_{\t_1,p_1})$ and
$(A_{\t_0,p_0},A_1)$.

We can also reformulate Proposition 1.11 using
duality. Assume $A_0\cap A_1$ dense in $A_0$ and in $A_1$
and also assume that there is a subspace $s \subset A_0
\cap A_1$ which is dense in $S_0$ with respect to $A_0$,
and in $S_1$ with respect to $A_1$. Then $(S_0, S_1)$ is
$K$-closed in $(A_0, A_1)$ iff $(S^\bot_0, S^\bot_1)$ is
$K$-closed in $(A^*_0, A^*_1)$.\medskip

\n {\bf Remark 1.13:} \ An even more abstract fact is
behind the preceding statement. Indeed, Proposition 1.10
\quad can be  viewed
as a consequence of the following statement: \ Let $X_1,
X_2$ be two closed subspaces of a $B$-space $X$. Let
$Q_i\colon \ X\to X/X_i$ be the quotient map $(i=1,2)$.
Then $Q_1(X_2)$ is closed iff $Q_2(X_1)$ is closed. This
can also be made more quantitative. Let us say that a
surjective operator is a $\lambda$-surjection if the image
of the open ball with center 0 and radius $\lambda$
contains the open unit ball with center 0. Now in the
above situation, if $Q_{1|X_2}$ is a $\lambda$-surjection
onto its image $Q_1(X_2)$, then $Q_{2|X_1}$ is a
$(\lambda+1)$-surjection onto its image $Q_2(X_1)$. (To
see the connection with Proposition 1.10, consider the
case $X=A_0\times A_1 , X_1=S_0 \times S_1 $ and  $X_2
=\{  (x,-x)\ \v \ x\in A_{0} \cap A_{1} \}$ .

\bigskip
\vfill\eject 
\magnification=\magstep1
\baselineskip=18pt
\def\ra{\rightarrow}

\def\n{\noindent}

\bigskip
\n{\bf 2. The non commutative case.}

Let $H$ be a  separable Hilbert space.	Let us denote by $C_p(H)$ or
simply by $C_p$ the Schatten ideal formed by all the compact
operators $T$ on $H$ such that $tr|T|^p<\infty$ and
equipped with the norm $\|T\|_{p}=(tr|T|^p)^{1/p}$.	Here
$1\leq p < \infty$.  If $p=\infty$, we denote by $C_\infty$
the space of all compact operators on $H$.  In the above,
we have taken $|T|$ defined as $(T^*T)^{1/2}$, but actually
this choice is unimportant here since (as is well known)
$tr(T^*T)^{p/2}=tr(TT^*)^{p/2}$, and hence
$\|T\|_{C_p} =\|T^*\|_{C_p}.$

Assume $H$ separable (possibly finite dimensional) and let $(e_n)$ be a
{\it fixed\/} orthonormal basis.  Let $E_k={\rm
span}(e_i ,i\le k)$ for all $k\geq 1$.	We will simply say
that a bounded operator $T:H\rightarrow H$ is triangular
if  $T(E_k)\subset E_k$ for all $k$. This definition can
be extended formally to the case when the indexing set
for the orthonormal  basis is any countable
totally ordered set in the place of the set of
all positive integers.

We will denote by $T_p(H)$ or simply by $T_p$ the subspace of $C_p(H)$
formed by all the triangular operators.  We will show the following
non-commutative version of P.~Jones' theorem proved in the preceding
section. The first point was proved recently, with $C_2$
and $T_2$ in the place of $C_1$ and $T_1$, in [KLW]. It
is this result from [KLW] which motivated the present
paper.

\proclaim Theorem 2.1.
\item{\rm (i)} There is a constant $K$	such that for any $x$ in $C_1$,
there is an operator $\widehat x$ in $T_1$ such that we have simultaneously
$$\eqalign{\|x-\widehat x\|_1&\leq K\ d_1(x,T_1)\cr
\|x-\widehat x\|_\infty&\leq K\ d_\infty(x,T_\infty)}$$
where we have denoted
$$d_p(x, T_p)=\inf\left\{\|x-y\|_p,y\in T_p\right\}.$$
\item{\rm (ii)} There is a constant $C$ such that for any $x$ in
$T_1+T_\infty$, we have
$$\forall t>0\quad K_t(x,T_1,T_\infty)\leq C\
K_t(x,C_1,C_\infty).$$ \item{\rm (iii)} If $0<\theta<1$ and
${1\over p}={{1-\theta}\over 1}+{\theta\over \infty}$, we
have $$(T_1,T_\infty)_{\theta p}=T_p$$
with equivalent norms (and similarly for the Lorentz space case).
\medskip

The proof of this theorem is entirely similar to the proof given in \S1, so
that we will only briefly review the main ingredients one by
one.

First we recall the following well known fact (cf.e.g. [GK])

\proclaim Lemma 2.2.  The orthogonal projection $P:C_2\ra
T_2$ is bounded {\it simultaneously\/} on $C_p$ for all
$1<p<\infty$.  Therefore, in particular {\rm (0.4)}
extends to the present non-commutative setting as
follows: If $1<p_0<p_1<\infty$ then there is a constant
$C$ such that for all $t>0$ we have $$\forall x\in
T_{p_0}+T_{p_1},\quad\forall t>0,\quad
K_t(x,T_{p_0},T_{p_1})\leq C\
K_t(x,C_{p_0},C_{p_1}).\leqno(2.1)$$ \medskip

The following fact is also well known.

\proclaim Lemma 2.3.  Let $1\leq p, q, r\leq \infty$ with
${1\over p}={1\over r}+{1\over q}$.  Assume (for
simplicity) that $H$ is finite dimensional.  Then every
invertible $x$ in $T_p$ can be written as  $x= ab$ with
$a\in T_r, b\in T_q$ and $\|a\|_r\|b\|_q=\|x\|_p$. \medskip

\noindent{\bf Proof:}  Note that either $p/q\leq 1/2$ or
$p/r\leq 1/2$. Assume $p/q \leq
1/2$.  Also, assume $\|x\|_{p}=1$.  By the Cholesky
factorization, we have $|x|^{2p/q}=b^*b$ for some $b$
triangular.  Moreover, $b$ is necessarily invertible and
$|b|^q=|x|^p$, so that $\|b\|_q=1$.  Let $\alpha =1-2p/q$.
Note $\alpha\geq 0$.  We have $x=U|x|$ with $U$ unitary.  Hence,
$x=U|x|^\alpha b^*b$.	Then, let $a=xb^{-1}=U|x|^\alpha
b^*$. Clearly, $a$ is triangular (since $a$ and $b^{-1}$
are so, and triangular operators form an algebra) and
moreover, $$\eqalign{aa^*&=U|x|^\alpha b^*b|x|^\alpha
U^*=U|x|^{2\left(1-{p\over q}\right)}U^*\cr &=U|x|^{2p\over
r}U^*.\cr}$$ Hence $(aa^*)^{1\over 2}=U|x|^{p\over r}U^*$
so that $\|a\|_r=1$.  This completes the proof if ${p\over
q}\leq {1\over 2}$. \medskip

In case we have instead ${p\over r}\leq {1\over 2}$, an entirely similar
argument works. We first write  $x=(xx^*)^{1\over 2}V$ then
$(xx^*)^{p\over r}=aa^*$ with $a$ triangular and
$b=a^{-1}x=a^*(xx^*)^{{1\over 2}-{p\over r}}V$.Then the rest
is the same. \  q.e.d. \medskip

It is now easy to complete the proof of the following non-commutative
version of Proposition 1.1, we skip the proof.
\medskip

\proclaim Proposition 2.4.  For all $1<p<q<\infty$ we have
$$T_p\subset (T_1,T_q)_{\theta,\infty}\leqno(2.2)$$
with norms bounded by some constant $K(p,q)$ where
$0<\theta<1$ satisfies ${1\over p}=1-\theta+{\theta\over
q}$.\medskip

To extend the squaring argument (0.1) (cf.the proof of
Proposition 1.2) , we need a non-commutative analogue of
the ``scaling'' that we used heavily for an analytic
function without zeros.  In the present setting, this is
somewhat easier.  Indeed, let us denote by
$\lambda(|x|)=(\lambda_n(|x|))_{n\geq 0}$ the sequence of
the eigenvalues of $|x|$ (arranged in non increasing
order, and repeated as usual according to their
multiplicity).	Observe that for all $\alpha >0$
$$\lambda_n(|x|^\alpha)=(\lambda_n(|x|))^\alpha.\leqno(2.3)\
\ $$

\proclaim Proposition 2.5.  If $0<p_0\leq p_1\leq \infty$
then for all $x$ in $C_{p_0}+C_{p_1}$ we have for all $t>0$
$$K_t\left(x,C_{p_0},C_{p_1})=K_t(|x|, C_{p_0},
C_{p_1})=K_t(\lambda(|x|),\ell_{p_0},\ell_{p_1}\right).$$
Moreover, a similar double identity holds for the $J$-functional.
\medskip

This is easy to check using the unitary invariance of the
spaces $C_p$ (for the first identity) and the existence of
a projection bounded on all $C_p$'s onto the elements
which are diagonal on the same basis as $|x|$ (for the
second one). In any case, we refer the reader to [PT] for
more details on such results. Using Proposition 2.5, it is
easy to extend Proposition 1.2 to the non-commutative case,
as follows.

 \proclaim Proposition 2.6.  For each $1<q<\infty$, there is
a constant $C^{'}$ (depending only on $q$) such that
$$\forall t>0\quad\forall x\in T_1+T_q\quad K_t(x,T_1,T_q)\leq
C^{'}\ K_t(x,C_1,C_q).$$\medskip

\noindent {\bf Proof:}	Assume w.l.o.g.~that $H$ is finite
dimensional and $x$ invertible.  Let $T$ be triangular such
that $$|x|=(x^*x)^{1\over 2}=b^*b=|b|^2.$$
We have $\lambda(|b|)=\lambda(|x|^{1\over 2})$.  Hence by Proposition 2.5
 and (2.3) we have
$$K_{t^{1\over 2}}(b,C_{2p_0},C_{2p_1})\leq
\left(2K_t(x,C_{p_0},C_{p_1})\right)^{1\over 2}.$$
Assume for simplicity $2K_t(x,C_{p_0},C_{p_1})<1$.  Then by (2.1), there
are $g_0,g_1$ triangular such that $b=g_0+g_1$ and $\|g_0\|_2+t^{1\over
2}\|g_1\|_{2q}<C$.  Note that $x=U|x|=Ub^*b=ab$ where $a=xb^{-1}$ is
triangular.  Since $a=Ub^*$ (and $\|x\|_p=\|x^*\|_p$ ) we
have obviously $$K_{t^{1\over
2}}(a,C_{2p_0},C_{2p_1})=K_{t^{1\over 2}}(b,C_{2p_0},
C_{2p_1}).$$ Hence by (2.1) again, $a=h_0+h_1$ with
$h_0,h_1$ triangular such that $\|h_0\|_2+t^{1\over
2}\|h_1\|_{2q}<C$.  Finally, $x=ab=(h_0+h_1)(g_0+g_1)$ can
be estimated as in the (commutative) proof of Proposition
1.2. q.e.d. \medskip

\noindent{\bf Remark 2.7:}  The analogue  of Proposition 1.3
is clearly valid in the case of $T_p$ with the same proof.
The same comment applies to Proposition 1.2$^*$ and
Proposition~1.3$^*$.  Moreover, Proposition 1.4 clearly
also extends to the non-commutative case, so that the proof
of Theorem 2.1 can be completed exactly as in \S1.

\vfill \eject

\def\i{\infty}
\magnification\magstep1
\baselineskip=18pt
\magnification\magstep1
\baselineskip=18pt

\def\i{\infty}
\def\l{\ell}
\def\t{\theta}
\def\la{\lambda}
\def\ra{\rightarrow}
\def\i{\infty}
\def\w{\widetilde}
          \def\u{\cal U}
\def\n{\noindent}
\def\e{\epsilon}

\def\ha{ \w{H}^{p_{0}}(C_{q_0}) }
\def\hb{ \w{H}^{p_{1}}(C_{q_1}) }
\def\LA{L_{p_{0}}(C_{q_0}) }
\def\LB{L_{p_{1}}(C_{q_1}) }
\n{\bf 3. The Banach space valued case.}\bigskip
We first remark that Jones' theorem remains valid for a
couple $(H^1(B),H^\i(B))$ for an arbitrary Banach space
$B$. Indeed, if $f$ is $H^1(B)$,  using an elementary
outer function argument, one can factor $f$ as $F\phi$,
 where $F$ is analytic, scalar valued and such that
$|F| = \|f\|_{B}$ , while $\phi$ is analytic $B$-valued
and
such that $\|\phi\|=1$ a. s. on {\bf T}. Moreover, if $f$
is in the unit ball of $\w{H}^1(B)$, then $f$ can be
approximated in  ${H}^1(B)$ by products $F\phi$ with
$F$ in the unit ball of $H^1$ and $\phi$ in the
unit ball of $\w{H}^\i (B)$. This reduces the
problem to the scalar case, since
it easy to verify that $K_t(f,H^1(B),H^\i(B))\le
K_t(F,H^1,H^\i) $ and
$K_t(F,L_1,L_\i)\le K_t(f,L_1(B),L_\i(B))$. Similarly, for
any $1\le p,q\le \i$, the couples $(H^p(B),H^q(B))$ and
$(\w{H}^p(B),\w{H}^q(B))$ are K-closed, by reduction to
the scalar case.

However, the general case of a compatible couple
$(A_0,A_1)$ of two different Banach spaces, is more
delicate. We refer the reader to [BX] for more information
and for a counterexample showing that Jones' theorem does
not extend in that degree of generality (cf. also
[X1],[X2]). Nevertheless, we show below that in a number
of nice cases, it does extend. There seems to be no
counterexample known at the time of this writing within
couples of Banach lattices.

Using the same method as in section 1, we can prove

\proclaim Theorem 3.1. Let $1\leq p_{0}, p_{1}, q_{0},
q_{1} \leq \infty$. Moreover, the space $\l_\i$ must be replaced by
$c_0$ wherever it appears.  Then the  couple
$(\w{H}^{p_{0}}(\ell_{q_0}),\w{H}^{p_{1}}(\ell_{q_1})) $
is K-closed with respect to
$(L_{p_{0}}(\ell_{q_0}),L_{p_{1}}(\ell_{q_1})) $.\medskip

\n {\bf Remark:} One can derive from Jones' proof the case
$p_{0}= q_{0}=1, p_{1}= q_{1}=\i$  , but
probably not the other cases. However, some other cases
can be derived from [B]. More precisely, Bourgain states
explicitly in [B] a theorem which in our terminology
means that the couple $(\w{H}^1(\l_1),\w{H}^1(\l_\i))$ is
K-closed. By a rather simple factorisation argument (such
as Theorem 2.7 in [HP]), one can show that a couple
$(\w{H}^p(A_0),\w{H}^p(A_1))$ is K-closed for some $1\le p\le \i$
iff it is K-closed for all $1\le p\le \i$. Therefore,
Bourgain's theorem does imply certain cases of Theorem
3.1. But actually, it is  interesting that one can
go conversely: in [P] we indeed do recover most of the
results of [B] by the methods of the
present paper.

We will denote simply by $g.h$ the pointwise and
coordinatewise product of two sequences $g=(g_n)$ and
$h=(h_n)$ of scalar analytic functions.
We first observe that in
the situation of Theorem 3.1, if $${1\over p} ={1\over
{2p_0}} +{1\over{2p_1}} \quad{\hbox {and}}\quad {1\over q}
={1\over {2q_0}} +{1\over{2q_1}},\leqno (3.1)$$ then the
unit ball of $H^p(\l_q)$ coincides with the set of all
products $g.h$ with $g$ and $h$  in the unit balls
respectively of  $H^{2p_{0}}(\ell_{2q_0})$ and
$H^{2p_{1}}(\ell_{2q_1})$. Indeed, this is easy to check
using outer functions.  Then, the squaring argument (0.1) suitably
generalized, becomes

\proclaim Lemma 3.2. In the same situation as in Theorem
3.1 (replacing  $\l_\i$  by
$c_0$ wherever it appears), if
$(\w{H}^{2p_{0}}(\ell_{2q_0}),\w{H}^{2p_{1}}(\ell_{2q_1}))$
is K-closed and if, with $p$ and $q$ as in (3.1) above,
$$\w{H}^p(\l_q) \subset
(\w{H}^{p_{0}}(\ell_{q_0}),\w{H}^{p_{1}}(\ell_{q_1}))_{1/2,\i}
,\leqno(3.2)$$ then
$(\w{H}^{p_{0}}(\ell_{q_0}),\w{H}^{p_{1}}(\ell_{q_1}))$ is
K-closed. \medskip

\n {\bf Proof:} The assumptions allow to reduce to the
case of finite dimensional $\l_p$ -spaces with constants
independent of the dimension. Then, the distinction
between ${H}^p(\l_q)$ and $\w{H}^p(\l_q)$ becomes
irrelevant and we can argue exactly as in Proposition
1.2. q.e.d.

  \n {\bf Proof of Theorem 3.1:} Let us record here the
preliminary observation that if we a priori know that
$(\w{H}^{2p_{0}}(\ell_{2q_0}),\w{H}^{2p_{1}}(\ell_{2q_1}))$
is K-closed, then (3.2)  holds iff the couple
$(\w{H}^{p_{0}}(\ell_{q_0}),\w{H}^{p_{1}}(\ell_{q_1}))$
is K-closed. Indeed, Lemma 3.2 gives the only if part,
and the converse is clear since $(\w{H}^{p_{}}(\ell_{q})$
is included into $L^p(\l_q)$, hence into the complex
interpolation space
$(L^{p_{0}}(\ell_{q_0}),L^{p_{1}}(\ell_{q_1}))_{1/2}$,
and a fortiori into
$(L^{p_{0}}(\ell_{q_0}),L^{p_{1}}(\ell_{q_1}))_{1/2,\i}$,
but, if we assume K-closedness, the latter space induces
on the subspace of analytic functions  a norm equivalent
to that of
$(\w{H}^{p_{0}}(\ell_{q_0}),\w{H}^{p_{1}}(\ell_{q_1}))_{1/2,\i}$,
which proves the converse part.

Using the observations
preceding Lemma 3.2, when all the indices are
finite it is easy to extend  Proposition
1.1 to the present setting with essentially the same
proof, more precisely, we have an inclusion $$
\w{H}^{p_\t}(\l_{q_\t}) \subset
(\w{H}^{p_{0}}(\ell_{q_0}),\w{H}^{p_{1}}(\ell_{q_1}))_{\t,\i}
,$$ where $ 1/{p_\t}=(1-\t)/p_0 +\t/p_1,\ \ 
1/{q_\t}=(1-\t)/q_0 +\t/q_1$, and  all the indices are
finite. Indeed, this can be checked using the well known
fact (apparently going as far back as [BB]) that the
Hilbert transform is bounded simultaneously on $L_p(\l_q)$
for all $1<p,q<\i,$ a fact which provides us with a
substitute for (0.4), and proves the preceding inclusion
when all indices are strictly between $1$ and $\i$ (this
can also  be seen, perhaps more easily,  using complex
interpolation). Then, the  factorisation argument as in
Proposition 1.1 yields the preceding inclusion assuming
only that all indices are finite. This extension and Lemma
3.2 give us Theorem 3.1 in case all the four indices $
p_{0}, p_{1}, q_{0}, q_{1}$ are finite. We now dualize. To
avoid irrelevant complications let us assume for the
moment that, everywhere until said otherwise, $\l_p$  is
the finite dimensional space ${\bf C} ^n$ equipped with
the $\l_p$-norm. Then, by Proposition 0.1 and dualisation,
we obtain Theorem 3.1 and (3.2) in case all the four
indices are more than $1$. At this stage, by the
preliminary remark recorded at the start of the proof, to
conclude it suffices to check that (3.2) holds in full
generality.   To do so, we note that if $q_0,q_1$ are
both finite and $p_0,p_1$ are both more than $1$, but
possibly infinite,  then we still have (3.2) and hence
K-closedness because we can apply the argument for
Proposition 1.1 to the second indices only. More
precisely, choosing $r$ large enough, we can write any
$f$ in the unit ball of  $\w{H}^{p}(\l_{q})$ as a product
$gh$ with $g$ in the unit ball of ${H}^{\i}(\l_{r})$ and
$h$ in the unit ball of $\w{H}^{p}(\l_{s})$ with $1/r
+1/s =1/q$, and this "translation " by $1/r$ reduces the
problem to the case of all indices more than $1$, which
has already been settled. By duality (or by a similar
argument applied to the first indices), if  $q_0,q_1$ are
both more than $1$ and $p_0,p_1$ both finite, we also
have K-closedness. Let us now check that all the other
cases follow. To describe the argument, it is convenient
to denote $x_0 =(1/p_0,1/q_0), x_1=(1/p_1,1/q_1) $ and to
view these two points as the extremities on a line
segment lying in the unit square of ${\bf R}^2$. Then
using the same argument as in Proposition 1.5 (or
invoking Wolff's theorem [W]) we can obtain (3.2) for the
segment $(x_0,x_1)$ everytime we know it for two
subsegments  $(x_0,y)$ and $(z,x_1)$ which intersect in a
non-empty open subsegment $(z,y)$. In this way, it is
then an entirely elementary matter to check all the
remaining cases, using the already settled ones.  This
concludes the proof in the case of finite dimensional
 $\l_p$-spaces, with constants independent of the
dimension. By a density argument, (note that the presence
of the tildas and the substitution of $c_0$ for $\l_\i$
allows the reduction to the finite dimensional case) it is
easy to deduce the general case from the finite dimensional
one.   q.e.d.\medskip

\proclaim Theorem 3.3. In the same situation as in
Theorem 3.1, the  couple
$(\w{H}^{p_{0}}(C_{q_0}),\w{H}^{p_{1}}(C_{q_1})) $ is K-closed
with respect to $(L_{p_{0}}(C_{q_0}),L_{p_{1}}(C_{q_1})) $.
\medskip

\n{\bf  Proof:} This is entirely analogous to the preceding
argument for Theorem 3.1, but of course we must use
suitable matrix-valued extensions of the
classical factorization theorems used in section 1.
Sarason' s paper [S] contains all that is needed here,
but actually, by density we need only prove the matrix
valued case, with constants independent of the size
of the matrix. In that case, if $H$ is finite dimensional,
it can be deduced from classical results of
Wiener-Masani-Helson-Lowdenslager (see [H]) that if
$1\le p_1,s_1,r_1 \le \i$ and $1/p_1=1/s_1+1/r_1$ , every
$f$ in the unit ball of $H^\i(C_{p_1}(H))$ can be written
as a product $f=gh$  with $g$ in the unit ball of
$H^\i(C_{s_1}(H))$ and $h$  in the unit ball of 
$H^\i(C_{r_1}(H))$. Here of course the product means the
pointwise product of operator valued functions. This can be
checked following the same argument as for Lemma 2.3 but
using the fact that any positive matrix valued function
$W$, such that (say) $W\ge \epsilon I$ for some $\epsilon
>0$, can be written as $b^* b$ for some  
invertible (actually outer) matrix valued analytic
function $b$ , cf. [H] . One can get rid of $\epsilon$ a
posteriori by a weak*-compactness
 limiting argument. Taking into account the remarks at the
beginning of this section, we find that, if $1\le
p_0,s_0,r_0 \le \i$ and $1/p_0=1/s_0+1/r_0$, any $f$ in
the unit ball of $H^{p_0}(C_{p_1}(H))$ can be written as a
product $f=gh$  with $g$ in the unit ball of
$H^{s_0}(C_{s_1}(H))$ and $h$  in the unit ball of 
$H^{r_0}(C_{r_1}(H))$. The proof of Theorem 3.3 can then
be completed easily following the same line of reasoning
as in section 2. Let us indicate here an
"economic" route for the inexperienced reader. We
assume $H$ 
finite dimensional, but all our constants will be
independent of its dimension. Let us explain more
technically what the  "squaring" argument becomes in the
non-commutative case. We will show that if 
$(\w{H}^{2p_{0}}(C_{2q_0}),\w{H}^{2p_{1}}(C_{2q_1}))$ is
K-closed and if (3.2) holds, then $(\ha,\hb)$ is K-closed.
To prove that consider $f$ in $\ha+\hb$ such that
$K_t(f,\LA,\LB) <1,$ this means there are $f_0\in \LA$ and
$f_1\in \LB$ such that $$f=f_0+f_1,\quad \|f_0\|_{\LA}
+ t\|f_1\|_{\LB} <  1. \leqno (3.3)$$
 Fix $\epsilon>0$. Let $F$ be an
analytic matrix valued function such that
$$F^{*}F=|f|+\epsilon I \leqno(3.4)$$ a.e. on ${\bf T}$,
and such that $z\ra F(z)^{-1}$ is analytic. This is
possible by choosing $F$ outer, cf.[H]. Let us denote by
$\u$ the set of all matrix valued (measurable) functions
$V$ such that $\|V(t)\| \le 1$ a.e. on ${\bf T}$. Note that
$|f|^{1/2}=VF$ for some $V\in \u.$ By polar decomposition,
$f=U |f| =U|f|^{1/2} V F$ for some $U\in \u$. Therefore,
$f=GF$ with $G= U|f|^{1/2} V$. But $G$ must be analytic
since $G=fF^{-1}$. Now we claim that for some $\delta >0$
which can be made arbitrarily small by letting $\e$ tend
to zero, we have
$$K_{t^{1/2}}(F,L_{2p_{0}}(C_{2q_0})+L_{2p_{1}}(C_{2q_1}))
< 2^{1/2}+ \delta ,\quad
K_{t^{1/2}}(G,L_{2p_{0}}(C_{2q_0})+L_{2p_{1}}(C_{2q_1}))
<2^{1/2} .\leqno (3.5)$$

Let us justify this. Going back
to our assumption on $f$, we have $|f|=U^{*}
f=U^{*}f_0 + U^{*}f_1$ . Using the diagonal projection
mentioned after Proposition 2.5 (more precisely, the
inequality $\|\sum a_{ii} e_i\otimes e_i\|_p \le
\|(a_{ij})\|_p $ valid for any $p\ge 1$ and any
orthonormal basis $(e_i)$),  we can project the last
decomposition and we obtain $|f|= g_0 + g_1$ where $g_0$
and $g_1$ are diagonal for the same basis as $|f|$ and by
(3.3) they satisfy $\|g_0\|_{\LA} + t \|g_1\|_{\LB} <1$.
Now since $|f| , g_0, g_1 $ all commute we can write
$|f|^{1/2} \le |g_0|^{1/2} + |g_1|^{1/2} $ from which it
follows exactly as in the commutative case that 
$$K_t(|f|^{1/2},L^{2p_{0}}(C_{2q_0}),L^{2p_{1}}(C_{2q_1}))
<2^{1/2} ,$$ which, recalling (3.4) and the value of $G$, 
obviously implies the above claim (3.5).  The rest of the
proof is then clear: since we assume K-closedness for the
doubled indices, we can write, for some constant $C$,
$$F=F_0 +F_1 \quad G=G_0 + G_1 $$ with
$$\|F_0\|_{\w{H}^{2p_{0}}(C_{2q_0})}+ t
\|F_1\|_{\w{H}^{2p_{1}}(C_{2q_1})} < C(2^{1/2}+\delta )$$
$$\|G_0\|_{\w{H}^{2p_{0}}(C_{2q_0})}+ t
\|G_1\|_{\w{H}^{2p_{1}}(C_{2q_1})} < C 2^{1/2}    .$$
Finally, we have $$F=GF= (G_0 F_0 + G_1 F_1)+ (G_0 F_1
+G_1 F_0 )$$ and we can conclude exactly as we did for
Proposition 1.1. This concludes the proof of the squaring
 argument. With this, it is now easy to complete the
proof exactly as in Theorem 3.1. q.e.d.

 \medskip

\proclaim  Corollary 3.4. Let $H$ be a 
separable Hilbert space. The couple $(H^1(C_1(H),
H^\i(B(H))$ is K-closed relative to  $(L_1(C_1(H)),
L_{\i}^w(B(H)))$, where we have denoted by
$L_{\i}^w(B(H))$ the space of all essentially bounded
weak*-measurable functions with values in $B(H)$
(if we view the $B(H)$-valued functions as matrix
valued, weak*-measurability  simply means here that all the
entries are measurable) .

\medskip

  \n {\bf Proof:} Consider  $f\in H^1(C_1(H))+
H^\i(B(H))$ with $K_t(f,L_1(C_1(H)),
L_\i(B(H))) <1$. We view $f$ as a doubly infinite matrix
valued function. Let $f_n$ be the function which has the
same entries as $f$ on the upper left $n\times n$ square
and zero elsewhere. By Theorem 3.3, for any $t>0$, we can
write $f_n =g_n +h_n$ with
$$ \|g_n\|_{H^1(C_1(H))}\le C \quad {\hbox {and} }\quad
\|h_n\|_{H^\i(B(H))} \le C/t .\leqno (3.6)$$
By  compactness, we can assume w.l.o.g. that $g_n$ and
$h_n$ converge in the weak operator topology, uniformly
on compact subsets of the unit disc $D$ to $g$ and $h$.
Clearly (3.6) remains valid in the limit and $f=g+h$, so
that we conclude $K_t(f,H^1(C_1(H)),
H^\i(B(H))) \le 2C $. By homogeneity, this completes the
proof. q.e.d.\medskip

By well known results on the interpolation of
$L_p$-spaces (cf.[BL] p.130, note 5.8.6, and [PT]) the
preceding results immediately  imply

\proclaim Corollary 3.5. If $1/p=1-\t,$ 
then $$(H^1(\l_1),H^\i(\l_\i))_{\t,p} = H^p(\l_p),$$
 and  $$(H^1(C_1(H)),
H^\i(B(H)))_{\t,p} =H^p(C_p(H))$$
for any separable Hilbert space $H$. Moreover, a similar
result holds for the $\w{H}^p$-spaces.

\vfill \eject
\medskip

\magnification\magstep1
\baselineskip=18pt

\def\i{\infty}
\def\l{\ell}
\def\t{\theta}
\def\la{\lambda}
\def\ra{\rightarrow}
\def\i{\infty}
\def\L{\Lambda}
\def\n{\noindent}
\def\w{\widetilde}
          \def\u{\cal U}
\def\n{\noindent}
\def\e{\epsilon}

\def\H{\bar{H}}
\def\v{\vert}
\def\V{\Vert}
\vfill\eject

 \n{\bf  4. Complex interpolation.}
\bigskip

In this section, we derive the complex version of Peter
Jones' theorem from the real one. Somehow, we
feel that the idea in the  proof of this deduction is of
some (theoretical) interest even for couples of $L_p$
spaces.

Let us denote by
$H^{p,q}_0$   the subspace of
$H^{p,q}$ formed by the boundary values of the analytic
functions vanishing at $0$. We will denote by  $\bar{H}^{p,q}_0$  the
subspace of $L_{p,q}$ formed by all the antianalytic
functions whose complex conjugates lie in $H^{p,q}_0$.
When $p=q$, as usual we denote these spaces by  $H^{p}_0
,\bar{H}^{p}_0$. We wish to prove that if  $1<p<\i$, and
$1/p = 1-\t$, then  $H^p =(H^1,H^\i)_\t .$
By standard
methods, it suffices to show that 
$$(L_1/\bar{H}^1_0,L_\i/\bar{H}^\i_0)_{1-\t} \subset H^q
\leqno(4.1)$$ where $1/q =\t$.

Let us denote by $dn$ the counting measure on the integers.
We will denote simply by ${\L}_{q,\i}$ the space
$L_{q,\i}(dm \otimes dn)$ and by $h^{q,\i}$
(resp. $\bar{h}^{q,\i}_0$ )the subspace
 formed by the elements $f(t,n)$ such that for each $n$,
the function $f(.,n)$ is in $H^{q,\i}$ (resp.
$\bar{H}^{q,\i}_0$).  By the method
of the preceding section,  it is easy to show that the
couple $(\bar{h}^{1,\i}_0,\bar{h}^{\i}_0)$ ( which, of
course is equivalent to the couple
$({h}^{1,\i},{h}^{\i})$) is K-closed with respect to
$({\L}_{1,\i},{\L}_{\i})$. Indeed, by Theorem 3.1 and 
reiteration (cf. Remark 1.12 ) we know that the couple 
$({h}^{2,\i},{h}^{\i})$ is K-closed, so that by the
squaring argument,  it suffices to check that
$${h}^{2,\i}\subset ({h}^{1,\i},{h}^{\i})_{1/2,\i} .$$ The
latter inclusion is clear since by Theorem 3.1, we have
${h}^{2,\i}= ({h}^{1},{h}^{\i})_{1/2,\i}$ and
${h}^{1}\subset  {h}^{1,\i}$. Thus, we have checked the
K-closedness of the couple
$(\bar{h}^{1,\i}_0,\bar{h}^{\i}_0)$ with respect to
$({\L}_{1,\i},{\L}_{\i})$.

In particular, by Proposition 0.1, this implies  we have a
simultaneous good lifting for the quotient spaces
$({\L}_{1,\i}/\bar{h}^{1,\i}_0,L_\i(\l_\i)/\bar{h}^{\i}_0)$.
From this, it is easy to deduce using the $J$-method that
the space
$({\L}_{1,\i}/\bar{h}^{1,\i}_0,L_\i(\l_\i)/\bar{h}^{\i}_0)_{1-\t,\i}$
can be identified with the space
${\L}_{q,\i}/\bar{h}^{q,\i}_0$ where $1/q= \t$.

By
interpolation and reiteration, the natural (i.e.
orthogonal) projection  is bounded from ${\L}_{q,\i}$ into
$h_{q,\i}$  if $1<q<\i$. Therefore,
${\L}_{q,\i}/\bar{h}^{q,\i}_0$ can simply be identified
with $h_{q,\i}$. The result of this discussion is the
following 
\proclaim Lemma 4.1. If  $1<q<\i$, and $1/q = \t$, then
there is a bounded natural inclusion
$$({\L}_{1,\i}/\bar{h}^{1,\i}_0,L_\i(\l_\i)/
\bar{h}^{\i}_0)_{1-\t,\i}\subset
h_{q,\i}.$$

We will now introduce a maping $J_q$ from $L_q/\H^q_0$
into  ${\L}_{q,\i}/\bar{h}^{q,\i}_0$ for all $1\le q\le
\i$, as follows. We start by defining a mapping $K_q :L_q
\ra \L_{q,\i}$ by setting
$$ \forall F\in L_q \quad K_q(F)(t,n) = n^{-1/q}
F(t).$$ For any  positive
real $x$, we denote 
 by $[x]$ the largest integer $n$ which is less than $x$.
Then, we have $$\sum_{n>0} t^q m\{n^{-1/q} \v F\v > t\} =
t^q\int [{\v F\v^q \over t^q} ]dm \le \int \v F\v^q dm. 
\leqno (4.2) $$ 
Moreover, the supremum of the left side over all
$t>0$ is equal to the right hand side (to check this,
simply let $t$ tend to $0$). Hence $\V K_q(F)\V =\V  F\V
$, so that $K_q$ has norm $1$. Note that $K_q$ obviously
maps $\H^{q}_0$ into $\bar{h}^{q,\i}_0$, therefore we may
define $J_q:L_q/\H^q_0\ra {\L}_{q,\i}/\bar{h}^{q,\i}_0$
as the mapping  canonically associated to $K_q$. For $f\in
L_q/\H^q_0$, if $F\in L_q$ is a representant of the
equivalence class of $f$,  then $(n^{-1/q} F)$ is a
representant of $J_q(f)$. The next result is a key
observation allowing us to deduce the complex 
interpolation theorem from the real one.
\proclaim Lemma 4.2. If\  $1<q<\i,\  1/q= \t $, the
operator $J_q$ defines a bounded mapping from
$(L_1/\H^1_0, L_\i/\H^\i_0)_{1-\t} $ into $h_{q,\i}$.

\n {\bf proof:} For any $z$ with $0<Re(z)<1$, let $J^z$ be
the operator defined exactly as $J_q$ but with $n^{z-1}$
in the place of $n^{-1/q}$. Then, by (4.2), if $Re(z)=0$,
$J^z$  is clearly a contraction from $L_1/\H^1_0$ into 
${\L}_{1,\i}/\bar{h}^{1,\i}_0$, and if $Re(z)=1$, it is a
contraction from $L_\i/\H^\i_0$ into 
$L_\i(\l_\i)/\bar{h}^{\i}_0$. Hence, by complex
interpolation (namely Stein's interpolation theorem for
analytic families of operators), $J_q=J^{1/q}$ is a
contraction from $(L_1/\H^1_0, L_\i/\H^\i_0)_{1-\t} $
into  $({\L}_{1,\i}/\bar{h}^{1,\i}_0,L_\i(\l_\i)/
\bar{h}^{\i}_0)_{1-\t}$, hence a fortiori (cf. e.g. [BL]
p.102, see also the following remark for a
technical precision)
  , it is bounded from $(L_1/\H^1_0,
L_\i/\H^\i_0)_{1-\t} $ into \ \ \ \
$({\L}_{1,\i}/\bar{h}^{1,\i}_0,L_\i(\l_\i)/
\bar{h}^{\i}_0)_{1-\t,\i}$, so that we can conclude the
proof by Lemma 4.1.   q.e.d.

\n{\bf Remark:} In the preceding argument, there is a
slight problem because ${\L}_{1,\i}$ is not normable, and
the complex interpolation method is usually developped in
the locally convex setting (see however [JJ]). This
difficulty can be circumvented easily. Indeed, let us
denote simply $Q_1 = L_\i(\l_\i)/
\bar{h}^{\i}_0$ . Let $B_0$ be the Banach  space of
all sequences of measurable functions $(x_n)$ such that
$\int \sup(n\v x_n\v) dm <\i$, equipped with the norm
$\V(x_n)\V =\int \sup(n\v x_n\v) dm$. We will denote by
$S_0$ the subspace of $B_0$ formed by the sequences $(x_n)$
such that $x_n \in \H^1_0$ for all $n$. Finally, we set
$Q_0=B_0/S_0$. We will use the observation   that
$B_0\subset {\L}_{1,\i}$ and this inclusion has norm one,
so that we also have  $Q_0 \subset
{\L}_{1,\i}/\bar{h}^{1,\i}_0$ with norm one. Then, the
preceding argument shows that $J_q$ is a contraction from
$(L_1/\H^1_0, L_\i/\H^\i_0)_{1-\t} $ into
$(Q_0,Q_1)_{1-\t}$, hence a fortiori it is bounded into 
$(Q_0,Q_1)_{1-\t,\i}$, and finally by the preceding
observation, into
$({\L}_{1,\i}/\bar{h}^{1,\i}_0,L_\i(\l_\i)/
\bar{h}^{\i}_0)_{1-\t,\i}$. In this manner, we have
managed to remain within Banach spaces.

We can now obtain the complex case of Peter Jones' 
 theorem as a consequence of the real case. 

\proclaim Theorem 4.3.
If  $1<p<\i$, and $1/p = 1-\t$, then 
$$H^p =(H^1,H^\i)_\t .$$

\n{\bf proof:} Let $p$ be the conjugate of $q$, so that
$1/p +1/q =1$. The inclusion $(H^1,H^\i)_\t \subset H^p$
is obvious. To prove the converse we dualize. Hence we
have to prove that  $$ (L_1/\H^1_0,
L_\i/\H^\i_0)_{1-\t} \subset H^q$$ where $1/q =\t$. By
Lemma 4.2, it suffices to show 
$$\forall x \in H^\i \quad \V x\V_{H^q} = \V ( 
n^{-1/q} x)\V_{h_{q,\i}} .$$
But this follows from the simple identity
$$\int \v x\v ^q dm = \sup_{t>0}\{ t^q \ \sum m(\v x\v >
t n^{1/q} )\} .\leqno (4.3)$$
This concludes the proof. (Note that (4.3) means that
$K_q$ is an isometric embedding of $L_q$ into 
$\L^{q,\i}$.)  q.e.d.

\proclaim Corollary 4.4. For any Banach space $B$, in the
 same situation as in Theorem 4.3, we
have 
$$\w{H}^p(B) =(\w{H}^1(B),\w{H}^\i (B))_\t ,\leqno (4.4)$$
and
$$H^p(B) =(H^1(B),H^\i(B))_\t .\leqno (4.5)$$

\n {\bf Proof:} The obvious inclusion 
$(\w{H}^1(B),\w{H}^\i (B))_\t \subset (L_1(B),L_\i
(B))_\t = L_p(B) $ implies $$(\w{H}^1(B),\w{H}^\i (B))_\t 
\subset \w{H}^p(B) .$$ Moreover, we recall that for any
$f$ in $H^p(B)$, and any $r<1$, the function $f_r$ defined
by $f_r(z) =f(rz)$ is clearly in $\w{H}^p(B)$, and $\V
f\V _{H^p(B)} = \sup_{0<r<1} \V f_r\V_{\w{H}^p(B)}$.
Using this, we
 obtain similarly  $$({H}^1(B),{H}^\i (B))_\t 
\subset {H}^p(B) .$$ To check the converse, by the
factorisation argument mentioned at the beginning of
section 3, we can write any $f$ in $H^p(B)$ as a product
$f=gh$, with $g\in H^p$ and $h\in H^\i (B)$. By theorem
4.3, $g$ belongs to $(H^1,H^\i)_\t$, hence, by
interpolation, since the multiplication by $h$ maps $H^1$
into ${H}^1(B)$ and  $H^\i$ into $H^\i (B)$, the
function $gh$ belongs to $({H}^1(B),{H}^\i (B))_\t$. This
completes the proof of (4.5). We leave  the rest of the 
proof to the reader.

The proof of
Theorem 4.3 extends with almost no change to the
non-commutative case (as was pointed out to me by
 QuanHua Xu), as follows

\proclaim Theorem 4.5.
If  $1<p<\i$, and $1/p = 1-\t$, then 
$$T_p =(T_1,T_\i)_\t \leqno (4.6)$$
and  $$(H^1(C_1(H)),
H^\i(B(H)))_{\t} =H^p(C_p(H))$$
for any separable Hilbert space $H$. Moreover, a similar
result holds for the $\w{H}^p$-spaces.

\n{\bf  Proof:} The argument is entirely similar to the
above. Let us indicate how (4.6) can be checked. First,
the inclusion $ (T_1,T_\i)_\t \subset T_p $ is obvious,
so that it suffices to prove the converse one. Let us
denote  $S_q(H)=C_q(H) \cap T_1^{\perp}$ and
$S_{q,\i}(H)=C_{q,\i}(H) \cap T_1^{\perp}$. To
abbreviate, we will sometimes write simply $S_q$
instead of $S_q(H)$. By duality, it suffices to show, in
analogy with (4.1) that $$(T_1/S_1,T_\i/S_\i)_{1-\t}
\subset T_q$$ where $1/q=\t= 1-1/p$. Let
$H\otimes \l_2$ be the Hilbert space which is the
Hilbertian tensor product of $H$ and $\l_2$.
 We define a
mapping $K_q$ from $C_q(H)$ into $C_{q,\i}(H\otimes \l_2)$
by letting
$$K_q(x)= \sum n^{-1/q}x\otimes \delta_n\otimes \delta_n $$
where we have denoted by $(\delta_n)$ the canonical basis of
$\l_2$. It is easy to check that $K_q$ is an isometric
embedding from $C_q(H)$ into $C_{q,\i}(H\otimes \l_2)$.
Let us denote $\hat{H}=H\otimes \l_2$. Since we assume
given an orthonormal basis $(e_n)$ in $H$, we can order
the basis $(e_n\otimes \delta_k)$  using the lexicographic
order, so that we can define as  usual the notion of a 
"triangular" operator on $\hat{H}$. Then obviously, 
$K_q$ maps $S_q$ into $S_{q,\i}(\hat{H})$, hence it
induces a mapping $J_q$ from $C_q/S_q$ into
$C_{q,\i}(\hat{H})/S_{q,\i}(\hat{H})$. We again denote by
$ J^z$ the same mapping but with $n^{1-z}$ in the place of
$n^{-1/q}$. By reasoning exactly as in Lemma 4.2, we can
show that $J_q$ is bounded from
$(C_1/S_1,C_\i/S_\i)_{1-\t} $ into $T_{q,\i}(\hat{H})$.
Let us denote by $(a_k(x))_{k\ge 0}$ the sequence of
singular numbers of an operator $ x$, (i.e. with the
notation of section 2, we have $a_k(x)=\la_k(\v x\v)$ ).
We observe that the sequence $(a_n(K_q(x)))$ coincides
with the  non-increasing rearrangement of the
collection $\{n^{-1/q} a_k(x) \v n\ge 1, k\ge 0 \}$. This
implies that $\V K_q(x)\V_{q,\i} =\V x\V_q$ . Hence, we
can argue exactly as for Theorem 4.3 above, and we obtain
(4.6). We leave  the rest of the  proof to the reader.

\medskip\vfill\eject
\noindent {\bf 
Remark:}\quad While we were completing the present paper,
we received a copy of a preprint by Paul M\"uller [M]
which contains a strikingly simple proof of Peter Jones'
theorem, or at least of Corollary 1.7 above, by an
extremely simple probabilistic stopping time argument. It
seems unlikely however that his idea will yield the 
non-commutative case.

\noindent {\bf Final Remark:}\quad The research for this
paper was motivated by a preprint  of Kaftal, Larson and
Weiss, where Proposition 1.2* and its non-commutative
analogue for nest algebras are proved for $p=2$ using an
operator algebraic method related to Arveson's distance
formula. The author is most grateful to David Larson for
showing him a copy of that paper and for stimulating
conversations.
\vskip24pt

\centerline {\bf References}\vskip6pt
\item {[B]} J.Bourgain, New Banach space properties of
the disc algebra and $H^\infty$, Acta Math. 152 (1984)
1-48. 

\item {[BB]} R.P.Boas and S.Bochner, On a theorem of
Marcel Riesz for Fourier series.J.London Math. Soc.
14(1939) 62-73.

\item {[BL]} J.Bergh and J.L\"ofstr\"om, Interpolation
spaces, An introduction, Springer Verlag 1976.

\item {[BS]} C.Bennett and R.Sharpley, Interpolation of
operators.Academic Press,1988.

\item {[BuD]} A.Bukhvalov and A.Danilevitch,Boundary
properties of analytic and harmonic functions with values
in a Banach space, Mat. Zametki 31 (1982) 203-214.
English translation: Mat. Notes 31 (1982) 104-110.

\item {[BX]} O.Blasco and Q.Xu, Interpolation between
vector valued Hardy spaces. Journal Funct. Anal.
 To appear.

\item {[E]} G.Edgar, Analytic martingale convergence,
Journal Funct. Anal. 69 (1986) 268-280.

\item {[G]} J.Garnett, Bounded Analytic Functions.
Academic Press 1981.

\item {[GK]} I.C.Gohberg and M.G.Krein, Introduction
to the theory of linear nonselfadjoint
operators, Transl. Math. Monogrphs, Amer. Math. Soc.
Providence, RI, 1969.

\item {[GR]} J.Garcia-Cuerva and J.L.Rubio de Francia.
Weighted norm inequalities and related topics. North
Holland, 1985.

\item {[H]} H.Helson, Lectures on invariant
subspaces.Academic Press, New-York 1964.

\item {[HP]} U. Haagerup and G. Pisier, Factorization of
analytic functions with values in non-commutative
$L_1$-spaces and applications. Canadian J. Math. 41
(1989) 882-906.

\item {[J1]} P.Jones, $L^\infty$ estimates for the
$\bar{\partial}$-problem in a half plane. Acta Math. 150
(1983) \nobreak{137-152}.

\item {[J2]} P.Jones, Interpolation between Hardy
spaces,  in: Conference on Harmonic Analysis in
honor of Antoni Zygmund, (edited by W.Beckner,
A.Calder\'on,R.Fefferman and P.Jones) 
 Wadsworth Inc.,1983, vol.2, p.437-451.

\item {[JJ]} S.Janson and P.Jones, Interpolation between
$H^p$ spaces: The complex method.
Journal Funct. Anal. 48 (1982) 58-80.

\item {[KLW]} V.Kaftal, D.Larson and G.Weiss,
Quasitriangular subalgebras of semifinite von Neumann
algebras are closed. Preprint. To appear.

\item {[M]} P.M\"uller, Holomorphic martingales and
interpolation of operators.

\item {[P]} G.Pisier, A simple proof of a theorem of Jean
Bourgain. Preprint.

\item {[PT]} A.Pietsch and H.Triebel,
Interpolationstheorie f\"ur Banachideale von
beschr\"ankten linearen operatoren, Studia Math. 31 (1968)
95-109.

\item {[S]} D.Sarason, Generalized interpolation in
$H^\infty$. Trans. Amer. Math. Soc.127 (1967) 179-203.

\item {[W]} T.Wolff. A note on interpolation spaces.
Springer Lecture Notes in Math. 908 (1982) 199-204.

\item {[X1]} Q.Xu, Applications du th\'eor\`eme de
factorisation pour des fonctions \`a valeurs op\'e
-rateurs.  Studia Math. 95 (1989) 273-292.

\item {[X2]} Q.Xu, Real interpolation of some Banach
lattices valued Hardy spaces.Preprint, Pub.\quad Irma,
Lille 1990 vol.20, 8.

\item {[X3]} Q.Xu, Elementary proofs of two theorems of
P.W.Jones on interpolation of Hardy spaces. Preprint,  
Pub.\quad Irma,
Lille, 1989.

\vskip12pt

Texas A. and M. University

College Station, TX 77843, U. S. A.

and

Universit\'e Paris 6
Equipe d'Analyse, Tour 46-0, 4\`eme \'etage,
4 Place Jussieu, 75230 

Paris Cedex 05, France

\end